\documentclass{article}
\usepackage{graphicx} 
\usepackage[utf8]{inputenc} 
\usepackage[T1]{fontenc}    
\usepackage{hyperref}       
\usepackage{url}            
\usepackage{booktabs}       
\usepackage{amsfonts}       
\usepackage{nicefrac}       
\usepackage{microtype}      
\usepackage{lipsum}		
\usepackage{graphicx}
\usepackage{natbib}
\usepackage{doi}
\usepackage[labelfont=bf]{caption}
\usepackage{bm}
\usepackage{graphicx}
\usepackage[space]{grffile}
\usepackage{latexsym}
\usepackage{textcomp}
\usepackage{longtable}
\usepackage{tabulary}
\usepackage{booktabs,array,multirow}
\usepackage{amsfonts,amsmath,amssymb}
\usepackage{natbib}
\usepackage{url}
\usepackage{hyperref}
\hypersetup{colorlinks=false,pdfborder={0 0 0}}
\usepackage{etoolbox}
\usepackage{mathtools}
\usepackage{tikz} 
\usepackage{tikz-3dplot}
\usetikzlibrary{arrows.meta}

\title{Early Warning Signals Can Vanish or Amplify: Dimensionality in Complex Systems}
\author{Susanne Ditlevsen$^{*1}$ and Peter Ditlevsen$^{*2}$\\{\small 1. Department of Mathematical Sciences, University of Copenhagen, Copenhagen, Denmark}\\{\small 2. Niels Bohr Institute, University of Copenhagen, Copenhagen, Denmark}\\
{\small * Emails SD: susanne@math.ku.dk, PD: pditlev@nbi.ku.dk}}

\date{}

\begin{document}

\maketitle

\begin{abstract}
Early warning signals (EWS), such as increasing variance and autocorrelation, are widely used to anticipate critical transitions associated with saddle-node bifurcations. However, real-world systems are often high-dimensional and multiscale, potentially altering the classical behavior of EWS. Here, we investigate how dimensionality, inertia, and red noise influence the detectability of EWS prior to tipping points. We show that when observations are not aligned with the critical direction, stable dynamics in orthogonal directions can mask EWS until very near the bifurcation. We further demonstrate that second-order dynamics with damping modify the autocorrelation structure and may either enhance or reduce signatures of critical slowing down. Finally, we study how colored noise influences EWS. Our results show that the presence and strength of EWS are not universal properties of tipping systems but depend critically on system geometry and stochastic forcing, implying that the absence of detectable EWS does not necessarily rule out an approaching critical transition.
\end{abstract}

\section{Introduction} 

Critical transitions and tipping phenomena are ubiquitous in complex dynamical systems, ranging from climate and ecology to neurons and applications in engineering. In many systems, abrupt transitions are associated with the disappearance of a stable equilibrium through a bifurcation, where the saddle-node and Hopf bifurcations are the generic co-dimension one bifurcations controlled by a single control parameter, $\lambda$. When the control parameter passes a critical value, $\lambda_c$, the former leads to an irreversible transition to an alternative state, while the latter results in oscillatory behavior. Here, we focus on the saddle-node bifurcation, which in colloquial terms is denoted a tipping point. A central question is whether such a transition can be anticipated from observations when the control parameter approaches the critical value slowly over time, i.e., prior to the tipping point. Over the past decades, the notion of early warning signals (EWS), in particular, loss of resilience reflected in increasing variance and critical slowing down reflected in increasing autocorrelation, has emerged as a promising framework for detecting forthcoming bifurcations \citep{dakos:2008,lenton:2008, scheffer:2009, dakos:2024reviewEWS}.

The concept of EWS originates from effectively one-dimensional stochastic differential equations (SDE's), which describe dynamical systems subject to noise, where the dynamics near the bifurcation are dominated by a single slow mode. Since the EWS, variance and autocorrelation, are inherently statistical concepts, it implicitly requires independent fast timescale fluctuations (noise) to monitor the stability of the dynamical steady state, which in this case is the statistical equilibrium state. 
In such a setting, the onset of critical slowing down follows directly from the vanishing stability of the equilibrium as the saddle-node bifurcation is approached. However, realistic physical and climate systems are inherently multidimensional and typically contain several interacting time scales. Although the underlying saddle-node structure, hysteresis, and tipping behavior may persist when embedded in a higher-dimensional system, the observable EWS need not behave as predicted for the low-dimensional system. 

To investigate the EWS, we assume the known or unknown complex dynamics governed by a set of non-autonomous and non-linear dynamical equations 
$$ \dot x = f(x,\lambda(t))+ \eta_t,$$ 
where a dot above a variable denotes derivative with respect to $t$. The state vector $x$ is potentially high dimensional, $\lambda(t)$ is a time varying control parameter (forcing), and $\eta_t$ is an additive stochastic noise representing, say, unresolved fast fluctuations (typically chaotic) variables $y$. 

Denoting $ x_i^*(\lambda)$ the attracting value of variable $x_i$ subject to dynamics with $\lambda(t)=\lambda$ constant, we may loosely assume that for a complex system to exhibit tipping: $\lim_{\epsilon\rightarrow0}(x^*_i (\lambda_c+\epsilon)- x^*_i(\lambda_c-\epsilon))\not=0$, where $\lambda_c$ is the critical forcing at the tipping point. The reason this can only be loosely defined is that, in the limit $\epsilon\rightarrow 0$, the system will experience noise-induced escape to the other state, thus, a simultaneous reduction of the noise intensity would be required. However, this loose definition will suffice in the following.  

There are several time scales involved in the system. Let $\tau_y$ and $\tau_x$ be (typical) autocorrelation times for a subset $y$ and a subset $x$ of the variables, respectively, and let $\tau_\lambda$ be the typical time scale for changes in $\lambda$,  
$$\tau_\lambda \sim \frac{\lambda_c-\lambda}{\dot \lambda}.$$ 
The implicit assumption for the validity and relevance of EWS is time scale separations such that $\tau_y\ll \tau_x \ll \tau_\lambda$. The first inequality ensures that a stochastic modeling of the variables $y$ is valid; the second inequality implies local quasi-stationarity. Consequently, moments and correlation functions such as $\langle x_i\rangle_{\lambda(t)}, \langle x_ix_j\rangle_{\lambda(t)}$ and $\langle x_i(t')x_i(t'+\tau)\rangle_{\lambda(t)}$, where $\langle\cdot \rangle$ denote expectation, are good approximations of the local equilibrium values. 

The time scale $\tau_x$ need not be unique, but may represent a collection of characteristic time scales, as is frequently the case in fast-slow systems. In the tipping system there will also be time scales associated with noise induced escape from one equilibrium state to another. These time scales should be long enough that EWS can be calculated prior to a noise-induced transition. This is a "race against the clock", similar to the problem of critical slowdown; the waiting time for a noise-induced tipping decreases as the critical threshold is approached.    

In what follows, it will be transparent what the relevant time scales are; we investigate how multidimensional embedding and time scale separation influence the detectability of EWS near saddle-node bifurcations. We do so by analyzing three simple cases of multi-dimensional systems that illustrate different aspects of the behavior of EWS and the detection of tipping points. In the following, we shall write the dynamical equations in differential form.  

We show that even when the tipping mechanism remains unchanged, the classical indicators of critical slowing down may emerge only very close to the bifurcation point. In systems with fast and slow interacting variables, the dominant observed dynamics can remain controlled by stable fast modes over most of the approach to tipping, thereby masking the weakening stability associated with the critical mode. As a consequence, EWS may appear substantially later than expected from corresponding one-dimensional models, reducing the relevance of EWS for predicting critical transitions.

\section{Saddlenode bifurcation embedded in a multidimensional system}

\begin{figure} 
\centering 
















\begin{tikzpicture}[scale=6]

\def\r{1}
\def\ang{30}

\coordinate (O) at (0,0);
\coordinate (X) at (1.1,0);
\coordinate (Y) at (0,1.05);

\coordinate (Z) at ({1.2*cos(\ang)},{1.2*sin(\ang)});

\draw[-{Latex[length=3mm]},thick,font=\large]
(O)--(X) node[below] {$X$};

\draw[-{Latex[length=3mm]},thick,font=\large]
(O)--(Y) node[left] {$Y$};

\draw[-{Latex[length=3mm]},ultra thick,blue,font=\large]
(O)--(Z) node[above right,blue] {$Z$};

\draw[dotted, thick]
({0.9*cos(\ang)},0)
--
({0.9*cos(\ang)},{0.9*sin(\ang)});

\draw[dotted, thick]
(0,{0.9*sin(\ang)})
--
({0.9*cos(\ang)},{0.9*sin(\ang)});

\draw[thick]
(0.22,0)
arc[start angle=0,end angle=\ang,radius=0.22];

\node[font=\Large] at
({0.27*cos(\ang/2)},
 {0.27*sin(\ang/2)})
 {$\theta$};

\draw[dashed, thick]
(0.9,0)
arc[start angle=0,end angle=90,radius=0.9];

\node[below,font=\large]
at ({0.83*cos(\ang)},0)
{$\cos(\theta)X_t$};

\node[left,font=\large]
at (0,{0.9*sin(\ang)})
{$\sin(\theta)Y_t$};

\node[below,font=\large]
at (0.92,0)
{$X_t$};

\node[left,font=\large]
at (0,0.9)
{$Y_t$};

\node[font=\large]
at (0.72,0.92)
{$Z_t=\cos(\theta)X_t+\sin(\theta)Y_t$};

\end{tikzpicture}


\caption{{\bf Illustration of the observation geometry in $\mathbf{R}^d$}. The critical saddle-node direction is denoted by $X$, while observations are made along the direction $Z$, which forms an angle $\theta$ with $X$. The direction $Y$ is orthogonal to $X$ within the plane spanned by $X$ and $Z$.
\label{fig:Zobservation}}
\end{figure}
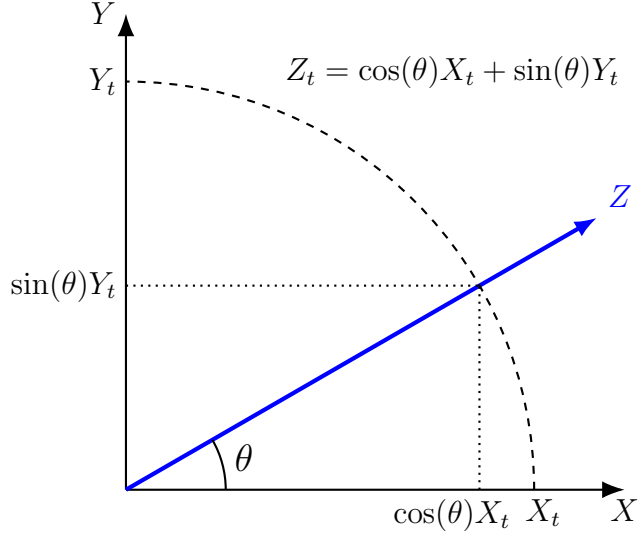

As the simplest extension to higher dimensions, we consider the situation where the system undergoes a saddle-node bifurcation along a one-dimensional critical direction in a $d$-dimensional state space \citep{abbott:2021,morr:2024SIAM}. Let $X$ denote a unit vector aligned with this critical direction. Suppose that observations are available only along a direction $Z$, whose alignment with the critical direction is unknown. Consider the two-dimensional subspace spanned by $X$ and $Z$. Within this subspace, define a coordinate system $(X,Y)$, where $Y$ is a unit vector orthogonal to $X$. Thus, $Y$ belongs to the orthogonal complement of $\mathrm{span}(X)$ and may be written as a linear combination of any basis of the $(d-1)$-dimensional subspace orthogonal to $X$. The observation direction can then be written as
\[
Z = \cos(\theta)\,X + \sin(\theta)\,Y,
\]
where $\theta \in [0,\pi/2]$ denotes the angle between the observation direction and the saddle-node direction, see Figure \ref{fig:Zobservation}. 

We assume that the $Y$ direction does not contain information about the approaching tipping point and is instead characterized by stationary fluctuations. Specifically, $Y$ is modeled as a zero-mean stationary Gaussian process. The saddlenode bifurcation is represented by its normal form. We have now reduced the $d$-dimensional state space flow to a canonical two-dimensional system in which the $X$ coordinate follows the normal form of a saddle-node bifurcation, whereas the $Y$ coordinate evolves as an independent zero-mean Ornstein-Uhlenbeck (OU) process.


The system is then the solution to the stochastic differential equation (SDE)

\begin{eqnarray}
dX_t &= -(X_t^2 + \lambda)\,dt + \sigma\, dW_t^{(1)}, \label{eq:SNmodel}\\
dY_t &= -\beta Y_t\,dt + \tau\, dW_t^{(2)}. \nonumber
\end{eqnarray}
Thus, the two coordinates are independent, $X_t$ has a saddlenode bifurcation for $\lambda \uparrow 0$ and 
$Y_t$ is a stationary OU process with mean 0 and variance $\tau^2/2\beta$. A linearization around the stable fixed point of the drift function of $X_t$ yields
\begin{equation}
X_t \approx -\alpha (\lambda) (X_t-\mu (\lambda) ) dt + \sigma dW_t
    \label{eq:OUapproxGeneral}
\end{equation}
with $\alpha (\lambda) = 2\sqrt{|\lambda|}$ and $\mu (\lambda) = \sqrt{|\lambda|}$. For fixed $\lambda \ll 0$ far from the critical value $\lambda_c = 0$, the variance of $X_t$ is $\sigma^2/2\alpha(\lambda)$. We thus obtain
\begin{eqnarray}
\label{eq:VarZ}
    \mbox{Var}(Z_t) &=& \cos^2 (\theta)\mbox{Var}(X_t)+ \sin^2 (\theta)\mbox{Var}(Y_t) = \frac{\sigma^2_X}{4\sqrt{|\lambda|}} + \frac{\sigma^2_Y}{2\beta},
\end{eqnarray}
where we have defined $\sigma^2_X = \cos^2(\theta)\sigma^2$ and $\sigma^2_Y = \sin^2(\theta)\tau^2$, and used that $X_t$ and $Y_t$ are independent. We therefore conclude that, except in the degenerate case $\theta=\pi/2$, the variance of the observed process diverges, 
\[ \mathrm{Var}(Z_t)\to\infty \qquad \text{as} \qquad \lambda \uparrow 0. \] 
Hence, the observation $Z_t$ exhibits the characteristic loss of resilience associated with the saddle-node bifurcation, even when observations are not aligned with the critical direction. However, the contribution from $Y_t$ may mask the EWS when the system is still far from the tipping point. In particular, for sufficiently negative values of $\lambda$, the constant term $\sigma^2_Y/2\beta$ dominates the variance, making the increase in $\mathrm{Var}(X_t)$ difficult to detect. As $\lambda$ approaches the bifurcation point, the variance of $X_t$ grows and eventually dominates the contribution from $Y_t$, rendering the constant term negligible.

Straightforward calculations also show that the $\Delta$-step autocorrelation of $Z_t$ is given by
\begin{eqnarray}
    \mbox{AC}_\Delta (Z_t) &=& \frac{\mbox{AC}_\Delta (X_t)\mbox{Var}(X_t)+ \mbox{AC}_\Delta (Y_t)\mbox{Var}(Y_t)}{\mbox{Var}(X_t)+\mbox{Var}(Y_t)} \nonumber \\
    &=& e^{-2\sqrt{|\lambda_t|} \Delta}\frac{\frac{ \sigma^2_X}{4\sqrt{|\lambda|}} }{\frac{\sigma^2_X}{4\sqrt{|\lambda|}} + \frac{\sigma^2_Y}{2\beta}} + e^{-\beta \Delta}\frac{  \frac{\sigma^2_Y}{2\beta}}{\frac{\sigma^2_X}{4\sqrt{|\lambda|}} + \frac{\sigma^2_Y}{2\beta}}.
    \label{eq:ACZ}
\end{eqnarray}


The first term in the final expression \eqref{eq:ACZ} is the product of two factors, both of which converge to one as $\lambda \uparrow 0$. The second term is the product of a constant and a factor that converges to zero. It follows that the autocorrelation converges to one as $\lambda \uparrow 0$. Therefore, the observed process $Z_t$ also exhibits the characteristic critical slowing down associated with the saddle-node bifurcation.

Note that this analysis neither requires $Y$ to be Gaussian nor assumes that the $(d-1)$-dimensional subspace orthogonal to the saddle-node direction can be represented by a single process. More generally, the orthogonal dynamics may be multivariate. The essential assumption is that these dynamics remain stationary with finite variance. Consequently, their contribution to the variance of the observation process $Z_t$ is bounded and independent of $\lambda$, and therefore does not alter the divergence induced by the critical saddle-node direction. 

\begin{figure}
    \centering
    \vspace{-1.5cm}
    \includegraphics[width = \linewidth]{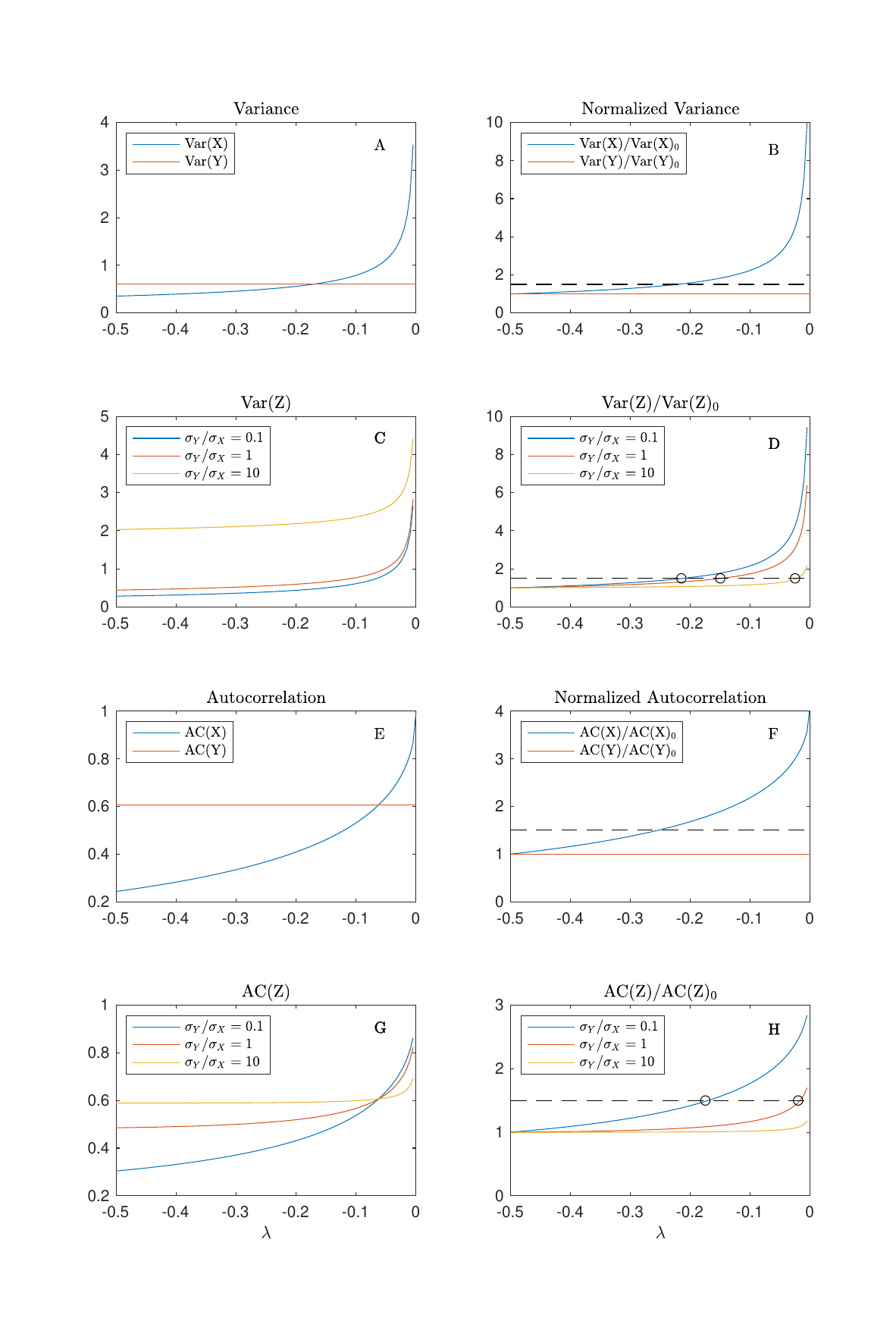}
    \vspace{-1.3cm}
    \caption{{\bf The observable $Z$ is a mixture of $X$ and $Y$}. \textbf{A:} Var$(X)$ increases as the bifurcation point $\lambda=0$ is approached, while Var$(Y)$ does not change with $\lambda$. \textbf{C:} Var$(Z)$ for $\sigma_Y=(0.1, 1, 10) \times \sigma_X$. \textbf{E} and \textbf{G}: Corresponding plots for the autocorrelation. \textbf{B}, \textbf{D}, \textbf{F}, \textbf{H:} The same for the EWS normalized with respect to $\lambda_0=-0.5$. The dashed lines show the confidence level at some significance, where the EWS has reached the level 50\% above the level for $\lambda_0$. The open circles in \textbf{D} and \textbf{H} indicate the up-crossings of the confidence level, showing that the smaller $\sigma^2_X$ is compared to $\sigma^2_Y$, the closer to the tipping point before an EWS can be detected.}
    \label{fig:MultidimensionalEmbedding}
\end{figure}

To illustrate the effect of masking the EWS by observing $Z$ rather than $X$, Figure \ref{fig:MultidimensionalEmbedding} shows the variance and autocorrelation as a function of the control parameter $\lambda$ for three cases where $\sigma^2_Y\ll \sigma^2_X, \sigma^2_Y=\sigma^2_X$ and $\sigma^2_Y\gg \sigma^2_X$. Only in the latter case, corresponding to conditions far from the tipping point, are the EWSs masked by the statistics of the process $Y$. As shown in Panels \textbf{D} and \textbf{F}, this substantially narrows the window for early warning, pushing detectable signals much closer to the tipping point. 

To investigate how the presence of $Y$ affects inference on the proximity to the tipping point, we consider the following scenario. Suppose that we treat the observed process as a direct proxy for $X_t$, thereby ignoring that the observations are, in fact, made along the direction $Z$. We then wish to estimate the time it takes for $\lambda$ to move from some baseline value $\lambda_0$ at time $t_0$ to a proportion $q, 0<q<1$ of the distance to the critical value, i.e., we wish to find $t$ such that
$$\lambda_{0} = q \lambda_t.$$
Ignoring the presence of $Y$, this corresponds to finding $t$ such that
$$\frac{1}{\mbox{Var}(Z_{t})^2} = q\frac{1}{\mbox{Var}(Z_{t_0})^2}$$
where we have used that 
$$|\lambda_t| = \left (\frac{\sigma^2}{4\mbox{Var}(X_t)} \right )^2$$
and we (wrongly) assume that $Z_t = X_t$. However, using \eqref{eq:VarZ} and after some algebra, we obtain
\begin{eqnarray}
    |\lambda_t| = q |\lambda_0| f(\lambda_0,\lambda_t)
\end{eqnarray}
with
\begin{eqnarray}
    f(\lambda_0,\lambda_t) &=& \frac{\frac{\sigma^4_X}{16}+\frac{\sigma^4_Y}{4\beta^2}|\lambda_t|+\frac{\sigma^2_X \sigma^2_Y}{4\beta}\sqrt{|\lambda_t|}}{\frac{\sigma^4_X}{16}+\frac{\sigma^4_Y}{4\beta^2}|\lambda_0|+\frac{\sigma^2_X \sigma^2_Y}{4\beta}\sqrt{|\lambda_0|}}
\end{eqnarray}
where the fraction $f(\lambda_0,\lambda_t)$ is less than $1$ whenever $\lambda_0 < \lambda_t<0$. Thus, at time $t$, the observed dynamics suggest that we have moved a fraction $q$ toward the tipping point, while the true fraction is $\tilde q = q f(\lambda_0,\lambda_t) < q$. Since smaller values of $q$ correspond to greater proximity to the tipping point, we are in fact closer to tipping than we believe, leading to an underestimation of the tipping risk due to $Y$.

\section{Relaxation times in tipping systems} 
\label{sec:SecondOrderSDE}

High dimensional complex systems with tipping may exhibit oscillations in variables that are either involved directly in the tipping or not. The former case may influence EWS for tipping. We shall not consider the case of Hopf-bifurcations but will examine a model exhibiting a saddle-node bifurcation including an acceleration term, which has damped oscillations up to some proximity to the tipping, where the stable point changes from being a spiral to a node. 
Consider the second-order equation
\begin{equation}
\ddot{x}_t + \eta \dot{x}_t = - \partial_x U(x_t, \lambda) + \sigma \xi(t), 
\label{eq:SecondOrder}
\end{equation}
where $U(x, \lambda)$ is a potential function, assumed to have a stable (possibly local) minimum $\mu (\lambda)$ and an unstable maximum for $\lambda < \lambda_c$, which coalesce in a saddlenode bifurcation for $\lambda = \lambda_c$. Parameter $\eta > 0$ is the damping, and $\xi(t)$ is white noise. We reformulate it as an SDE for variables $X_t$ and $V_t = \dot{X}_t$:
\begin{equation} \label{eq:SecondOrderSDE}
    \begin{aligned}
d X_t &= V_t d t, \\
d V_t &= \left (-\eta V_t - \partial_x U(X_t, \lambda) \right) d t + \sigma d W_t,
\end{aligned}
\end{equation}
where $W_t$ denotes a standard Wiener process. The invariant density is given by 
\begin{equation*}
    \pi(x,v; \lambda) = C \exp\left(-\frac{2 \eta}{\sigma^2}U(x, \lambda)\right) \exp\left(-\frac{\eta}{\sigma^2} v^2\right), \label{eq:XVinv}
\end{equation*}
where $C$ is the normalizing constant. The marginal invariant probability of $V_t$ is thus Gaussian with zero mean and variance $\sigma^2/(2\eta)$. The marginal invariant probability of $X_t$ is shaped by the potential $U(x, \lambda)$. 
At steady state and for fixed $\lambda$, the position $X_t$ and velocity $V_t$ are independent. This is reflected by the decomposition of the joint density $\pi(x,v; \lambda)$ into $\pi(x; \lambda) \pi(v)$.

\begin{figure}[ht]
    \centering
    \includegraphics[width = 10cm]{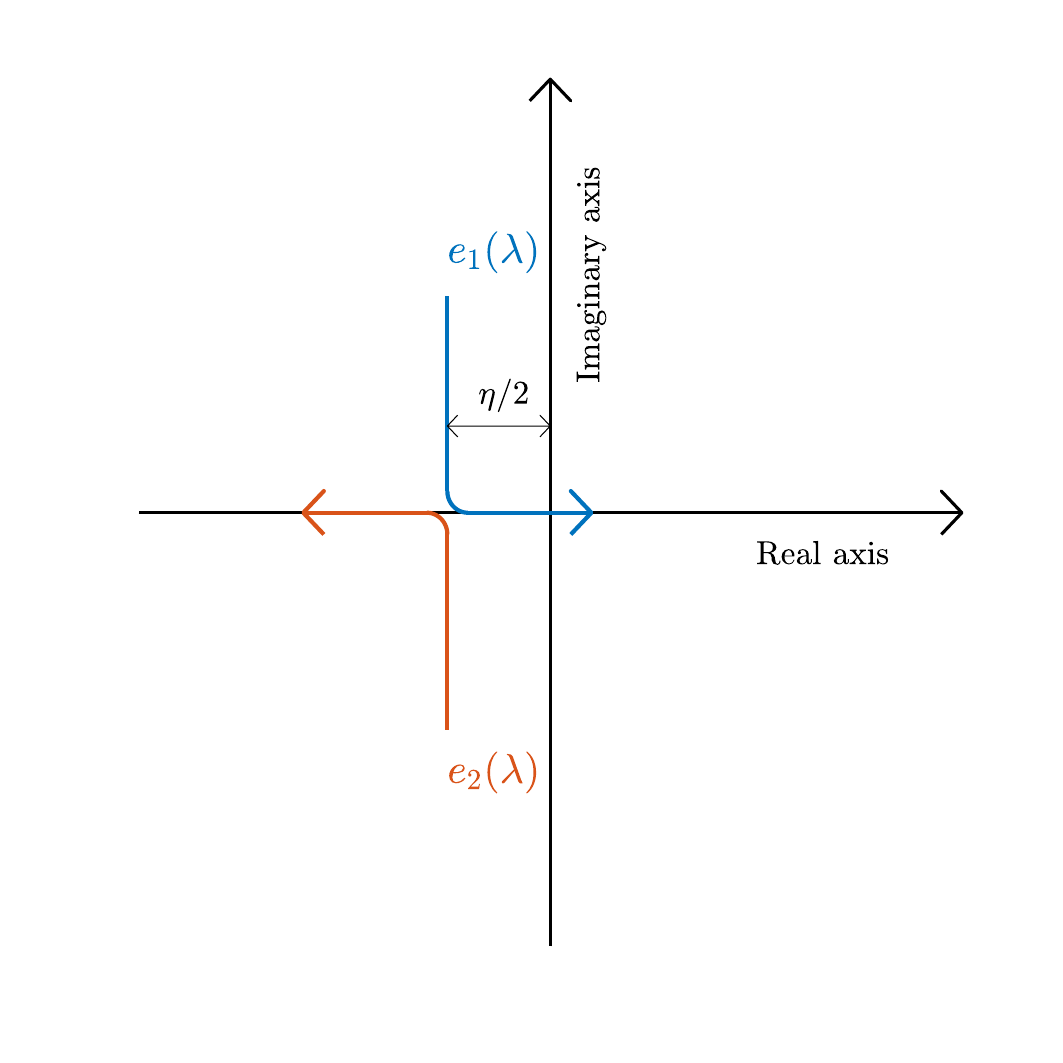}
  
    \caption{{\bf Stability of the second order saddlenode bifurcation model.} The colored curves show how the eigenvalues of the linearized dynamics (eigenvalues of the Jacobian) change as the control parameter $\lambda$ increase. For $\lambda<\lambda^*$, where $\lambda^*$ is defined by $\eta^2=4\alpha(\lambda^*)$, the equilibrium state is a stable focus with damped oscillations, while for $\lambda^* <\lambda<\lambda_c$ it is a stable fixed point. Tipping occurs for $\lambda=\lambda_c$ when the eigenvalue $e_1$ crosses the imaginary axis.}
    \label{fig:fig4}
\end{figure}

To find an approximation of the variance and autocorrelation of process $Y_t = (X_t, V_t)^\top$ for fixed $\lambda$, we Taylor expand the drift around the stable point $\mu (\lambda)$, to obtain an approximating OU process, from which we obtain closed expressions of approximations of variance and autocorrelation. This yields 

\begin{equation} \label{eq:OUapprox}
d Y_t \approx A(\lambda) (Y_t - b(\lambda)) d t + \Sigma d W_t,
\end{equation}
where $b(\lambda) = (\mu (\lambda), 0)^\top$ and $\Sigma = (0, \sigma)^\top$. The Taylor expansion yields
\begin{equation} \label{eq:A}
A(\lambda) = \begin{bmatrix}
    0 & 1\\
    - \alpha (\lambda) & -\eta\\
\end{bmatrix}
\end{equation}
where $\alpha (\lambda) = \partial_{xx} U(\mu (\lambda), \lambda)$. 
The eigenvalues of $A$ are
$$e_{1,2} = -\frac{\eta}{2} \pm i \frac{\sqrt{\eta^2 - 4\alpha (\lambda)}}{2},$$
which are complex conjugates for $D(\lambda) := \eta^2 - 4\alpha (\lambda) <0$ and real for $D(\lambda) \geq 0$.
At the bifurcation we have $\alpha(\lambda)\rightarrow 0$ as $\lambda \rightarrow \lambda_c$ corresponding to $e_1\rightarrow 0$ as illustrated in Figure \ref{fig:fig4}.

The covariance matrix of \eqref{eq:OUapprox} for a step size $\Delta$ is given by
\begin{eqnarray*}
    \Omega_\Delta &=& \int_0^\Delta
    e^{ A(\Delta - s)} \Sigma \Sigma^\top e^{ A^\top(\Delta - s) } ds \\
&=&\frac{\sigma^2}{(e_2-e_1)^2}\begin{bmatrix}
 \Omega_{\Delta,11}&\Omega_{\Delta,12}\\
 \Omega_{\Delta,21}& \Omega_{\Delta,22}\\
\end{bmatrix}
\end{eqnarray*}
with
\begin{eqnarray*}
    \Omega_{\Delta,11} &=& \frac{1}{2e_1}(e^{2e_1 \Delta}-1)+\frac{1}{2e_2}(e^{2e_2 \Delta}-1) +\frac{2}{\eta}(e^{-\eta\Delta}-1)\\
    \Omega_{\Delta,12} \, = \, \Omega_{\Delta,21} &=&\frac{1}{2}(e^{2e_1 \Delta}-1)+\frac{1}{2}(e^{2e_2 \Delta}-1) -(e^{-\eta\Delta}-1)\\
    \Omega_{\Delta,22} &=& \frac{e_1}{2}(e^{2e_1 \Delta}-1)+\frac{e_2}{2}(e^{2e_2 \Delta}-1) +\frac{2e_1e_2}{\eta}(e^{-\eta\Delta}-1)
\end{eqnarray*}
where we have used that
$$e^{At} = \frac{1}{(e_2-e_1)} \begin{bmatrix}
    1 & 1\\
    e_1 & e_2
\end{bmatrix}
\begin{bmatrix}
e^{e_1 t} & 0 \\
0 & e^{e_2 t}
\end{bmatrix}
\begin{bmatrix}
    e_2 & -1\\
    -e_1 & 1
\end{bmatrix}.$$
Letting $\Delta \rightarrow \infty$ we obtain the stationary variance:
\begin{eqnarray*}
    \Omega (\lambda) &=& \sigma^2
    \begin{bmatrix}
    \frac{1}{2\eta \alpha (\lambda)}
     & 0\\
    0 & \frac{1}{2\eta} 
    \end{bmatrix}.
\end{eqnarray*}
Thus, only the variance of $X_t$ diverges as the bifurcation is approached (bifurcation at $\alpha = 0$), while the variance of $V_t$ remains constant. The variance in the one-dimensional model with no second order term is given by
$$\mbox{Var} (X_t^{[1D]}) = \frac{\sigma^2}{2\alpha(\lambda)} = \frac{1}{\eta} \mbox{Var} (X_t^{[2D]}),$$
if $\sigma$ and potential are the same in the two models. For ease of notation, we suppress the dependence on $\lambda$ on $A,D$ and $\alpha$. The $\Delta$-lag autocovariance is given by $e^{A\Delta}\Omega (\lambda)$, but since $\Omega(\lambda)$ is diagonal, the $\Delta$-lag autocorrelation is directly given by
$$e^{A\Delta} = \frac{1}{(e_2-e_1)} 
\begin{bmatrix}
e_2 e^{e_1 \Delta}-e_1 e^{e_2 \Delta} & e^{e_2 \Delta}-e^{e_1 \Delta}\\
e_1 e_2 (e^{e_1 \Delta}-e^{e_2 \Delta}) & e_2 e^{e_2 \Delta}-e_1 e^{e_1 \Delta}
\end{bmatrix}.$$
For complex eigenvalues (i.e., $D<0$), the one-lag correlation then equals
\begin{align*}
&e^{A\Delta} = \\
&e^{-\frac{\eta \Delta}{2}} 
\begin{bmatrix}
\cos \left (  \frac{\sqrt{|D|}}{2} \Delta\right ) + \frac{\eta}{\sqrt{|D|}} \sin \left (  \frac{\sqrt{|D|}}{2} \Delta\right )  & \frac{2}{\sqrt{|D|}} \sin \left (  \frac{\sqrt{|D|}}{2} \Delta\right )\\
-\frac{2 \alpha}{\sqrt{|D|}} \sin \left (  \frac{\sqrt{|D|}}{2} \Delta\right ) & \cos \left (  \frac{\sqrt{|D|}}{2} \Delta\right ) - \frac{\eta}{\sqrt{|D|}} \sin \left (  \frac{\sqrt{|D|}}{2} \Delta\right )
\end{bmatrix}.
\end{align*}
Close to the bifurcation point ($\alpha (\lambda) < \eta^2/4$), eigenvalues will always be real, and
the one-lag autocorrelation equals
\begin{align*}
&e^{A\Delta} =\\
&e^{-\frac{\eta \Delta}{2} }
\begin{bmatrix}
\cosh \left (  \frac{\sqrt{D}}{2} \Delta\right ) + \frac{\eta}{\sqrt{D}}\sinh \left (  \frac{\sqrt{D}}{2} \Delta\right ) & -\frac{2}{\sqrt{D}} \sinh \left ( \frac{\sqrt{D}}{2} \Delta\right )\\
\frac{2\alpha}{\sqrt{D}} \sinh \left ( \frac{\sqrt{D}}{2} \Delta\right ) & \cosh \left (  \frac{\sqrt{D}}{2} \Delta\right )-\frac{\eta}{\sqrt{D}}\sinh \left (  \frac{\sqrt{D}}{2} \Delta\right )
\end{bmatrix}.
\end{align*}
The first term in the matrix corresponding to $X_t$ converges to $\cosh (\eta \Delta/2)+\sinh (\eta \Delta/2)= e^{\eta \Delta/2}$ as $\alpha \downarrow 0$, and thus, the autocorrelation converges to 1 as $\lambda$ approaches $\lambda_c$. Furthermore, for $\eta^2 \gg 4\alpha$, the autocorrelation is approximately equal to 1. Thus, for large $\eta$, it will be difficult to detect an increase in autocorrelation for increasing $\lambda$, since the autocorrelation is dominated by $\eta$, and close to 1. The autocorrelation is always larger for the 2D than for the 1D model and it increases faster when approaching the tipping point, and it is thus easier to detect critical slowing down when there is friction in the system. 

The last term in the matrix corresponding to $V_t$ converges to $\cosh (\eta \Delta/2)-\sinh (\eta \Delta/2)= e^{-\eta \Delta/2}$ as $\alpha \downarrow 0$, and thus, the autocorrelation converges to $e^{-\eta \Delta}$, and will not show critical slowing down.

As an example, we take $U(x,\lambda)$ to be the double-well potential, relevant for key tipping elements in climate \citep{Stommel:1961,cessi:1994}, which results in Kramers oscillator \citep{ArnoldImkeller2000}, which is the stochastic Duffing oscillator \citep{duffing1918erzwungene}. A quasi-likelihood parameter estimator for this model was proposed in \cite{pilipovic2025} and used to model paleoclimate data from the Greenland ice core. It has
\begin{equation}
U(x) = \frac{x^4}{4}-\frac{x^2}{2} + \lambda x.  \label{eq:Duffing}
\end{equation}
The SDE becomes
\begin{equation} \label{eq:KramersSDE}
    \begin{aligned}
d X_t &= V_t d t, \\
d V_t &=  \left (-\eta V_t - X_t^3 + X_t - \lambda \right) d t + \sigma d W_t.
\end{aligned}
\end{equation}
For $|\lambda|<\lambda_c = 2/\sqrt{27}$ there are three fixed points of the drift in \eqref{eq:KramersSDE}, two stable $(\mu_+,0)$, $(\mu_-,0)$ and one unstable $(\mu_0,0)$, where $\mu_+(\lambda) = 3\lambda_c\cos(\varphi (\lambda)/3)$, $\mu_0 (\lambda) = 3\lambda_c\cos((\varphi (\lambda)-2\pi)/3)$ and $\mu_-(\lambda) = 3\lambda_c\cos((\varphi (\lambda)-4\pi)/3)$, and $\varphi (\lambda) = \arccos (-\lambda/\lambda_c)$. 

\begin{figure}[ht]
    \centering
    \includegraphics[width = \linewidth]{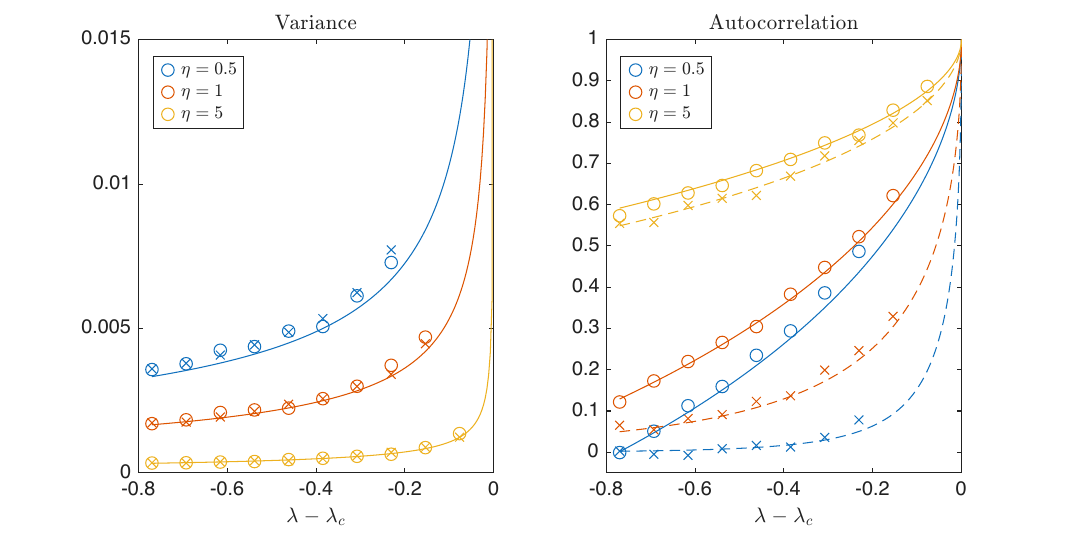}
    \caption{{\bf Saddlenode model with acceleration.} Variance and autocorrelation of $X_t$ for $\eta=0.5, 1, 5$, circles are obtained from simulations. Crosses are simulations of the first order model with no acceleration term. The corresponding curves are the analytic expressions (dashed for the first order model). The autocorrelation is always larger for the second order model compared to the first order model. There are no simulations for $\lambda$ too close to $\lambda_0$ since in this range the system will experience noise induced transition to the other minimum of the potential in a time scale faster than the time window necessary for detecting an EWS. Notice that this range is larger for smaller values of the damping $\eta$. 
    }
    \label{fig:acceleration}
\end{figure}

Setting $\mu = (\mu_+,0)$, we have $\alpha(\lambda) = 3 \mu_+^2-1$ and 
$$\mbox{Var} (X_t^{[DW,2D]}) = \frac{\sigma^2}{2\eta(3 \mu_+^2-1)}.$$
Variance and autocorrelation of $X_t$ generated from models with and without acceleration are illustrated in Figure \ref{fig:acceleration}. We set $\sigma_{[DW,1D]} = \sigma_{[DW,2D]}/\eta $ and $\alpha_{[DW,1D]} = \alpha_{[DW,2D]}/\eta $, where $_{[DW,1D]}$ denote the one-dimensional double-well potential model without acceleration and likewise for the two-dimensional model with acceleration, such that variances will be the same and autocorrelations comparable in the two models. The model without acceleration is simulated by setting the acceleration term $\ddot x_t$ to zero in \eqref{eq:SecondOrder}.



\section{Saddlenode model driven by colored noise}

The last example we will consider is colored noise  \citep{morr:2024}. We modify eq. \eqref{eq:OUapproxGeneral} by letting the noise term be an OU process, 
\begin{eqnarray}
\label{eq:OUapproxRedNoiseX}
X_t &\approx& -\alpha (\lambda) (X_t-\mu (\lambda) ) dt + Y_t dt \\
dY_t &=& -\beta Y_t dt + \sigma dW_t
    \label{eq:OUapproxRedNoiseY}
\end{eqnarray}
with $\alpha (\lambda) \rightarrow 0$ when $\lambda \rightarrow \lambda_c$. The Jacobian is 
\begin{equation} \label{eq:A2}
A(\lambda) = \begin{bmatrix}
    - \alpha (\lambda) & 1\\
    0 & -\beta\\
\end{bmatrix}
\end{equation}
with real eigenvalues $\alpha$ and $\beta$ with corresponding eigenvectors $(1,0)$ and $(1/(\alpha - \beta), 1)$. Similar calculations as in Section \ref{sec:SecondOrderSDE} yield
\begin{eqnarray*}
    \Omega (\lambda) &=& \frac{\sigma^2}{2 \beta }
    \begin{bmatrix}
    \frac{1}{\alpha (\lambda) (\alpha (\lambda)+\beta)}
     &  \frac{1}{(\alpha (\lambda)+\beta)} \\
    \frac{1}{(\alpha (\lambda)+\beta)} & 1 
    \end{bmatrix},
\end{eqnarray*}
\begin{eqnarray*}
    e^{A(\lambda)\Delta}  &=& 
    \begin{bmatrix}
     e^{-\alpha (\lambda)\Delta}
     &  \frac{1}{\alpha (\lambda)-\beta} (e^{-\beta\Delta}-e^{-\alpha\Delta}) \\
    \frac{1}{2 \beta (\alpha (\lambda)+\beta)} & 
    e^{-\beta\Delta} 
    \end{bmatrix}
\end{eqnarray*}
and
\begin{eqnarray*}
    e^{A(\lambda)\Delta}\Omega (\lambda) &=& 
    \begin{bmatrix}
    COV_{11} &  COV_{12}
      \\
    COV_{21} & COV_{22}     
    \end{bmatrix}
\end{eqnarray*}
with
\begin{eqnarray*}
    COV_{11} &=& \sigma^2\frac{\alpha (\lambda)e^{-\beta\Delta} - \beta e^{-\alpha(\lambda)\Delta}}{2 \alpha (\lambda)\beta (\alpha (\lambda)+\beta)(\alpha (\lambda)-\beta)}, \\
    COV_{12} &=& \sigma^2\frac{(\alpha (\lambda)+\beta)e^{-\beta\Delta}-2\beta e^{-\alpha(\lambda)\Delta}}{(\alpha (\lambda)+\beta)(\alpha (\lambda)-\beta)}, \\
    COV_{21} &=&  \sigma^2\frac{e^{-\beta\Delta}}{2 \beta (\alpha (\lambda)+\beta)},\\
    COV_{22} &=& \sigma^2\frac{e^{-\beta\Delta}}{2\beta}.
\end{eqnarray*}
The variance and autocorrelation of $X_t$ are thus 
$$\mbox{Var}(X_t) = \frac{\sigma^2}{2 \alpha (\lambda)\beta (\alpha (\lambda)+\beta)} \rightarrow\infty$$
and 
\begin{eqnarray*}
\mbox{AC}_\Delta(X_t) = \frac{COV_{11}}{\mbox{Var}(X_t)} &=&  \frac{\alpha (\lambda)e^{-\beta\Delta} - \beta e^{-\alpha(\lambda)\Delta}}{\alpha (\lambda)-\beta} \rightarrow 1
\end{eqnarray*}
when $\lambda \rightarrow 0$. Thus, EWS are still present, but are modulated by the relation between the timescales given by $1/\beta$ and the rate of change of $\alpha$ towards the tipping at 0. 

\begin{figure}[ht]
    \centering
    \includegraphics[width = \linewidth]{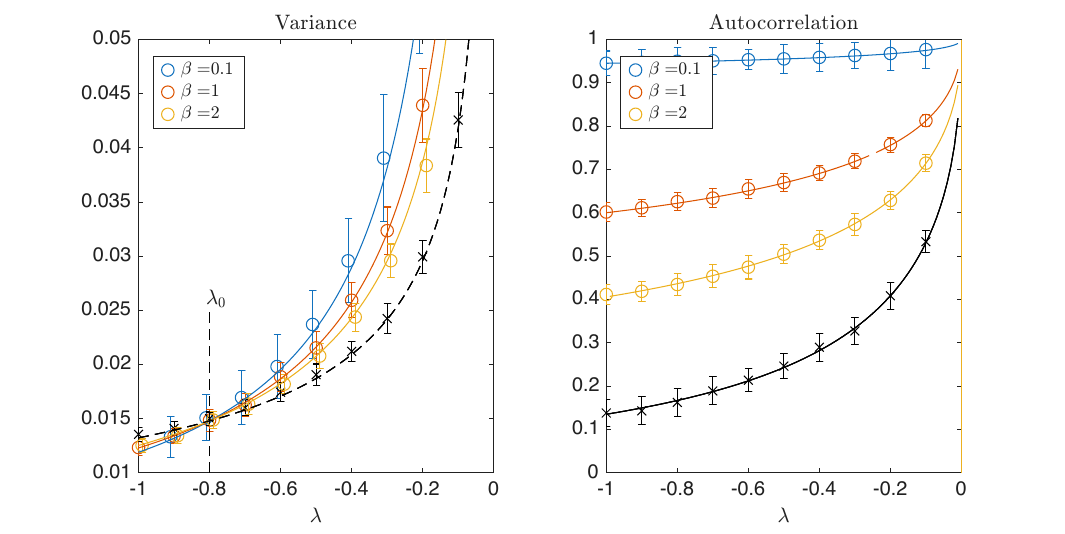}
    \caption{{\bf Saddlenode model driven by red noise.} Variance and autocorrelation of $X_t$ in \eqref{eq:OUapproxRedNoiseX} for $\beta=0.1, 1, 2$, circles are obtained from simulations. Crosses are simulations of the model driven directly by white noise. The corresponding curves are the analytic expressions (black dashed for the first order model). Parameters $\sigma$ are chosen such that the variances are the same for $\lambda_0 = -0.8$. The autocorrelation is always larger for the model driven by red noise compared to the model driven by white noise. The error bars provide one standard deviation from 100 simulations.
    }
    \label{fig:rednoise}
\end{figure}

Figure \ref{fig:rednoise} shows variance and autocorrelation for three values of $\beta$ and compared to the model driven by white noise, eq. \eqref{eq:OUapproxGeneral}. Parameters $\sigma$ are chosen such that the variances are the same for $\lambda_0 = -0.8$, i.e., $\sigma^2 (\beta) = \beta (\alpha(\lambda_0) - \beta)$ for the red noise models, and $\sigma^2 = \sigma^2(\beta)/\beta (\alpha(\lambda_0) - \beta)$ for the white noise model. 


The variances increase more rapidly for the red-noise models, with the rate of increase becoming larger as $\beta$ decreases. However, larger variances are also associated with stronger autocorrelation, implying that larger window sizes are required to detect EWS. Consequently, despite the slower increase in variance under white noise, EWS may still be easier to detect in this case because the increase can be estimated more reliably and is therefore statistically easier to identify.

\section{Discussion}

EWS are often presented as generic indicators of approaching critical transitions, with increasing variance and autocorrelation expected to arise as a consequence of critical slowing down. Our results demonstrate that this classical picture is not universally valid. Although the underlying saddle-node bifurcation remains unchanged, the appearance, strength, and timing of EWS depend strongly on how the system is observed and on the nature of the stochastic forcing acting upon it.

First, we showed that in multidimensional systems, observations need not be aligned with the critical direction. Stable dynamics in orthogonal directions can substantially mask the growth in variance and autocorrelation associated with the approaching bifurcation, causing detectable EWS to emerge only very close to the tipping point. This effect reduces the practical prediction horizon and may lead to an underestimation of tipping risk when observations are interpreted through low-dimensional models.

Second, incorporating inertia through second-order dynamics alters the expected behavior of autocorrelation. While critical slowing down remains present in the position variable, damping can either enhance or obscure increases in autocorrelation depending on the timescale structure of the system. 

Third, colored noise modifies the quantitative expression of EWS while preserving the asymptotic tendency toward increasing variance and autocorrelation.

Our analysis highlights an important limitation of conventional EWS theory when applied to complex multiscale systems. The results suggest that the absence of detectable EWS far from the tipping point does not necessarily imply the absence of an underlying critical transition, but may instead reflect the geometric, dynamical and timescale  structure of a higher-dimensional system. Consequently, the absence of clear EWS should not be interpreted as evidence that a system is far from a critical transition. Future work should focus on inference methods that exploit full dynamical information rather than relying solely on summary statistics such as variance and autocorrelation. Such approaches may offer a more robust framework for identifying loss of stability in complex, noisy, and high-dimensional systems before conventional EWS become detectable.

These findings furthermore suggest that improved statistical methodologies are needed to uncover the finer dynamical structure underlying multiscale systems and to detect weakening stability before classical summary-based EWS become visible. In particular, likelihood-based inference methods for nonlinear stochastic processes, such as those proposed by \cite{pilipovic2024} for general nonlinear stochastic systems, and applied to saddle-node tipping point models in \cite{Ditlevsenx2}, provide a promising framework for resolving dynamical changes that may remain hidden to conventional early warning statistics until very close to the bifurcation point.

\bibliographystyle{plainnat}
\bibliography{references}       

\end{document}